\documentclass[a4paper,10pt,onecolumn]{article}

\usepackage{a4wide,times,bm}
\usepackage{subeqnarray}
\usepackage{enumerate,graphicx,amsmath,amssymb,latexsym,psfrag,epsfig}
\usepackage{theorem,fancyhdr,fancybox}
\usepackage{calrsfs}
\usepackage[normalem]{ulem}
\usepackage{lscape}

\newtheorem{theorem}{Theorem}
\newtheorem{lemma}{Lemma}

\newtheorem{definition}{Definition}

{\theorembodyfont{\upshape}}
{\theorembodyfont{\upshape}}
{\theorembodyfont{\upshape}}

\begin{document}

\title{A Round-Robin Protocol for Distributed Estimation with $H_\infty$
  Consensus
\thanks{This work was supported under Australian Research Council's
  Discovery Projects funding scheme (project number DP120102152), and by
  Israel Science Foundation (grant No 754/10). Part of this work was
  carried out during the first author's visit to the Australian National
  University.}
\thanks{A version of this paper has been presented at the 52nd IEEE CDC,
  Florence, Italy.}}
\author{
 V.~Ugrinovskii\thanks{School of Engineering
   and Information
   Technology,  UNSW Canberra at the
 Australian Defence Force Academy, Canberra ACT 2600, Australia,
 Email: v.ugrinovskii@gmail.com, FAX: +61 2 6268 8443, Tel: +61 2 6268 8219.}
\quad and \quad
E. Fridman\thanks{School of Electrical Engineering, Tel-Aviv
  University, Tel-Aviv 69978, Israel, Email: emilia@eng.tau.ac.il.}}

\maketitle\thispagestyle{empty}

\begin{abstract}
 The paper considers a distributed robust estimation problem over a network
 with directed topology involving continuous time observers. While
 measurements are available to the observers continuously, the nodes interact
 according to a Round-Robin rule, at discrete time instances. The results
 of the paper are sufficient conditions which guarantee a suboptimal $H_\infty$
 level of consensus between observers with sampled
 interconnections.  
\end{abstract}

\section{Introduction}

The problem of distributed estimation is one of very active topics in the
modern control theory and signal processing literature. Interest in this
problem is motivated by a growing number of applications where a decision
about the observed process must be made simultaneously by spatially
distributed sensors, each taking partial measurements of the process.

When the process and measurements are subject to noise and disturbance, 
robustness aspects of the problem come into prominence. In the past several
years, a number of results have been presented in the literature which
develop the $H_\infty$ control and estimation theory for distributed
systems subject to uncertain perturbations; e.g.,
see~\cite{DLG-2011,LDCH-2010,LJL-2008,NF-2009,U6,LaU1,U7}. In particular,
methodologies of distributed sampled-data $H_\infty$ filtering have been
considered~\cite{SWL-2011}. That reference emphasized several 
aspects of realistic sensor networks, among them coupling between
sensor nodes through the information communicated between neighbouring
sensor nodes and the sampled nature of that coupling, which is dictated by
the digital communication technology. 

In this paper we address some of the above challenges by developing a
Round-Robin protocol for a network of distributed estimators. 
The Round-Robin protocol is a commonly used protocol for information
transmission in networked control systems. From a hybrid systems
perspective this protocol has been studied in details
in~\cite{NT-2004,HTWN-2010}. More recently, it has been considered in the
context of time-delay systems in~\cite{LFH-2012}, where an analysis of 
exponential stability and $L_2$ properties of networked control systems
with Round-Robin scheduling was presented using a delay switching system
modeling. In this paper, we further develop this technique in the context
of robust distributed estimation with intermittent communication between
sensing nodes.  

The objective of this paper is to develop a sufficient condition to enable
synthesis of filter and interconnection gains for a network of
distributed observers, whose aim is to track dynamics of a linear
uncertain plant. 

The first contribution of this paper is a version of the Round-Robin protocol
of~\cite{LFH-2012} to be used with the distributed estimation schemes
proposed~\cite{U6,LaU1,U7}. We show that instead of continuously exchanging
information 
(the type of networks considered in those references), 
the node observers can achieve the relative $H_\infty$ consensus objective
by exchanging information at certain sampling times, by polling one
neighbour at a time. 
Our second contribution demonstrates that the Round-Robin design
of~\cite{LFH-2012} can be applied to derive a network of
non-switching observers.
 
Our main result is a sufficient condition, expressed in the form of Linear
Matrix Inequalities (LMIs), from which filter and
interconnection gains for each node estimator can be computed, to ensure
the network of sampled data observers under consideration converges to the
trajectory of the observed plant. As
in~\cite{U6,LaU1,U7}, our methodology relies on
certain vector dissipativity properties of the
large-scale system comprised of the observers' error
dynamics~\cite{HCN-2004}. However, 
different form these references, to establish these vector dissipativity properties, we employ a novel class of generalized supply rates which reflect the
sampled-data nature of interconnections between observers. The general idea
behind introducing such generalized supply rates can be traced
to~\cite{LCD-2004} (also, see~\cite{LaU1}), but our proposal here makes use
of special properties of sampled signals. In the limit, when the maximum
sampling period approaches zero, these generalized supply rates vanish, and one
recovers the vector dissipativity properties of error dynamics established
in~\cite{U6}.  

The paper is organized as follows. The problem formulation, along with the
graph theory preliminaries is presented in Section~\ref{Formulation}. The
main results of the paper are given in Section~\ref{main}.
Section~\ref{Concl} concludes the paper.   

\paragraph*{Notation} Throughout the paper, $\mathbf{R}^n$ denotes a real
Euclidean $n$-dimensional vector space, with the norm $\|x\|\triangleq
(x'x)^{1/2}$; here the symbol $'$ denotes the transpose of a matrix or a
vector. $L_2[0,\infty)$ will denote the Lebesgue space of
$\mathbf{R}^n$-valued vector-functions $z(\cdot)$, defined on the time
interval $[0,\infty)$, with the norm $\|z\|_2\triangleq
\left(\int_0^\infty\|z(t)\|^2dt\right)^{1/2}$ and the inner product
$\int_0^\infty z_1(t)'z_2(t)dt$. $\otimes$ is the Kronecker product of
matrices, $\mathbf{1}_n\in \mathbf{R}^n$ is the column-vector of ones.

\section{The problem formulation}\label{Formulation}

\subsection{Graph theory}

Consider a filter network with
$N$ nodes and a directed graph topology $\mathcal{G} = (\mathcal{V},\mathcal{E})$;
$\mathcal{V}=\{1,2,\ldots,N\}$, $\mathcal{E}\subset \mathcal{V}\times
\mathcal{V}$ are the set of vertices and the set of edges, respectively.
The notation $(j,i)$ will denote the edge
of the graph originating at node $j$ and ending at node $i$. In accordance
with a common convention~\cite{OM-2004}, we consider graphs without
self-loops, i.e., $(i,i)\not\in\mathbf{E}$. However, each node is assumed
to have complete information about its filter and measurements.

For each $i\in \mathcal{V}$, we denote $\mathcal{V}_i=\{j:(j,i)\in
\mathcal{E}\}$ to be the \emph{ordered} set of nodes supplying information
to node $i$, i.e, the neighbourhood of $i$.  Without loss of generality,
suppose the elements of $\mathcal{V}_i$ are ordered in the ascending
order. The cardinality of $\mathcal{V}_i$, known as the in-degree of node
$i$, is denoted $p_i$; i.e., $p_i$ is equal to the number of incoming edges
for node $i$. Also, the out-degree of node $i$ (.e., the number of outgoing
edges) is denoted $q_i$.

Without loss of generality the graph $\mathcal{G}$ will be assumed to be
weakly connected, 
\cite[Proposition~1]{U6}.

In the sequel,  a shift permutation operator defined on elements of the set
$\mathcal{V}_i$ will be used:
\begin{equation}
  \label{permu}
  \Pi\{j_1,\ldots,j_{p_i-1},j_{p_i}\}=\{j_{p_i},j_1,\ldots,j_{p_i-1}\}.
\end{equation}
Furthermore, $\Pi^k(\mathcal{V}_i)$ will denote the set obtained from
$\mathcal{V}_i$ using $k$ consecutive shift permutations (\ref{permu}).
In regard to this set, the following notation will be used throughout the
paper unless stated otherwise: for $\nu\in \{1,\ldots, p_i\}$,
$j_\nu$ is the $\nu$-th element in the ordered set
$\Pi^k(\mathcal{V}_i)$. Conversely, for $j\in \Pi^k(\mathcal{V}_i)$,
$\nu_j^{k,i}\in \{1,\ldots, p_i\}$ is the index of element $j$ in the
permutation $\Pi^k(\mathcal{V}_i)$. We will omit the superscript $^{k,i}$
if this does not lead to ambiguity.

\subsection{Distributed estimation with $H_\infty$
  consensus}\label{Estimation}

Consider a plant described  by the equation
\begin{equation}
  \label{eq:plant}
  \dot x=Ax+B_2w(t), \qquad x(t)=x_0~\forall t\le 0. 
\end{equation}
Here $x\in\mathbf{R}^n$ is the state of the plant, and
$w\in\mathbf{R}^{m_w}$ is a disturbance. 
We also assume $w(t)\in L_2[0,\infty)$, so that the $L_2$-integrable solution
of (\ref{eq:plant}) exists on any finite interval
$[0,T]$~\cite[p.125]{CP-1977}.

The distributed filtering problem under consideration is to
estimate the state of the system (\ref{eq:plant})
using a network of filters connected according to the graph
$\mathcal{G}$. Each node takes measurements
\begin{equation}\label{U6.yi}
y_i(t)=C_ix(t)+D_{2i}w(t)+\bar D_{2i}v_i(t);
\end{equation}
$v_i\in\mathbf{R}^{m_v}$ is a measurement disturbance. 

The measurements are processed by a network of observers connected over the
graph $\mathcal{G}$. The key assumption in this paper is to allow the
observers make use of their local 
 measurements continuously, however they can only interact with each other
 at discrete time instances $t_k$, $k=0,1\ldots$, with $t_0=0$. For
 simplicity, we assume that this schedule of updates is known to all
 participants in the network, and therefore all nodes
 exchange information at the same time instance $t_k$. However, at every
 time instance $t_k$ only one neighbour in the set $\mathcal{V}_i$ is polled by
 each node $i$, according to the `Round-Robin' rule. Formally, this leads
 us to define the following observer protocol: For $t\in[t_k,t_{k+1})$,
 $k=0,1,\ldots$,
 \begin{eqnarray}
    \dot{\hat x}_i&=&A\hat x_i(t) + L_i(y_i(t)-C_i\hat x_i(t)) \nonumber \\
         &&+K_i\sum_{j\in \Pi^k(\mathcal{V}_i)} H_i(\hat x_j(t_{k-\nu_j+1})-\hat
           x_i(t_{k-\nu_j+1}),
  \label{UP7.C.d}
\end{eqnarray}
where $ \hat{x}_i(t)$ is the estimate of the plant
state $x(t)$ calculated at node $i$,
the matrices $L_i$,
$K_i$ are parameters of the filters to be determined, and $H_i$ is a given
matrix. The observers are initiated with zero initial condition, $\hat
x_i(t)=0$ for all $t\le 0$.

From now on, we will omit the time variable when a signal is
considered at time $t$, and will write, for example, $\hat x_i$ for $\hat
x_i(t)$. 

The last term in (\ref{UP7.C.d}) reflects the desire of each node observer to
update its estimate of the plant using feedback from the neighbours in
its neighbourhood, according to the consensus estimation
paradigm~\cite{Olfati-Saber-2007,U6}. However, unlike these references,
under the Round-Robin protocol, only one neighbour is polled
at each time $t_k$ to provide a `neighbour feedback', and this sample is
stored and used by the observer until time $t_{k+p_i}$. The feature of the
Round-Robin protocol is to poll the neighbours one at a time, in a cyclic
manner. Formally, this can be described 
by first applying the shift permutation operator $\Pi$ to the neighbourhood
set at every time instance $t_k$, and then selecting the first element from the
resulting permutation $\Pi^k(\mathcal{V}_i)$ for feedback.

Let $e_i=x-\hat x_i$ be the local estimation error at node $i$. This error
satisfies the equation:
 \begin{eqnarray}
    \dot{e}_i&=&(A - L_iC_i)e_i+(B-L_iD_i)\xi_i \nonumber \\
      &+&K_iH_i \sum_{\j\in \Pi^k(\mathcal{V}_i)}
      (e_{j_\nu}(t_{k-\nu_j+1})-e_i(t_{k-\nu_j+1})),
\label{e.1} \\
&& \qquad \qquad\qquad\quad t\in [t_k, t_{k+1}), \ k=0,1,\ldots. \nonumber
\end{eqnarray}
Here we used the notation $\xi_i$ to represent the perturbation vector
$[w'~v_i']'$, and the matrices $B$, $D_i$ are defined as follows $B=[B_2~0]$,
$D_i=[D_{2i}~\bar D_{2i}]$. The initial conditions for (\ref{e.1}) are 
$e_i(t)=x_0$ $\forall t\le 0$. In particular in (\ref{e.1}),
$e_{j_\nu}(t_{k-\nu_j+1})-e_i(t_{k-\nu_j+1})=0$ for $t_k-\nu_j+1<0$.  

Since the error dynamics (\ref{e.1}) are governed by $L_2$ integrable
disturbance signals $\xi_i$, we can only expect the node observers to converge
in $L_2$ sense. To quantify transient consensus performance of the observer
network (\ref{UP7.C.d}) under disturbances, consider the cost of
disagreement between the observers caused by a particular vector of disturbance
signals $\xi(\cdot)=[\xi_1(\cdot)'~\ldots~\xi_N(\cdot)']'$,
\begin{eqnarray}
  \label{Psi.L2}
  J(\xi)&=&\frac{1}{N} \int_0^\infty \sum_{i=1}^N \sum_{j\in
  \Pi^k(\mathcal{V}_i)}\|\hat x_j(t)-\hat x_i(t)\|^2 dt \nonumber  \\
&=& \frac{1}{N} \int_0^\infty \sum_{i=1}^N \sum_{j\in
  \Pi^k(\mathcal{V}_i)}\|e_j(t)-e_i(t)\|^2dt,
\end{eqnarray}
where $k$ is a time-dependent index, $k=0,1,\ldots$, defined so that for every
$t\in[0,\infty)$, $t_k\le t<t_{k+1}$. The functional (\ref{Psi.L2}) was
originally introduced in~\cite{U6} as a measure of consensus performance of
a corresponding continuous-time observer network. It is worth noting that
for each $t$, $\sum_{j\in \Pi^k(\mathcal{V}_i)}\|\hat x_j(t)-\hat
x_i(t)\|^2$ is independent of the order in which node $i$ polls its
neighbours, so that
\[
\sum_{j\in \Pi^k(\mathcal{V}_i)}\|\hat x_j(t)-\hat
x_i(t)\|^2=\sum_{j\in \mathcal{V}_i}\|\hat x_j(t)-\hat
x_i(t)\|^2.
\]
Therefore, the inner summation in
(\ref{Psi.L2}) can be replaced with summation over the neighbourhood set
$\mathcal{V}_i$. This observation leads to the same expression for $J(\xi)$
as in the case of continuous-time networks~\cite{U6},   
\begin{equation}
  \label{Psi.L2.1}
  J(\xi)=\frac{1}{N} \int_0^\infty \sum_{i=1}^N\left
    [(p_i+q_i)\|e_i(s)\|^2-2e_i'\sum_{j\in \mathcal{V}_i}e_j(s) \right] ds.
\end{equation}

The following distributed estimation problem is a version of the
distributed $H_\infty$ consensus-based estimation problem originally
introduced in~\cite{U6,LaU1}, modified to include the Round-Robin protocol
(\ref{UP7.C.d}).

\begin{definition}\label{Def1}
The distributed estimation problem under consideration is to
determine a collection of observer gains
$L_i$ and interconnection coupling gains
$K_i,~i=1,\ldots,N$, for the filters
(\ref{UP7.C.d}) which ensure that the following conditions are
satisfied:
  \begin{enumerate}[(i)]
  \item
In the absence of uncertainty, the interconnection of unperturbed
systems (\ref{e.1}) must be exponentially stable.

\item
The filter must ensure a specified level of transient consensus performance, as
follows 
\begin{eqnarray}\label{objective.i.1}
&&\sup_{x_0, \xi\neq 0}
\frac{J(\xi)}{\|x_0\|^2_P+\frac{1}{N}\|\xi\|_2^2} \le
\gamma^2. 
\end{eqnarray}
Here, $\|x_0\|^2_P=x_0'Px_0$, $P=P'>0$ is a matrix to be determined,
and $\gamma>0$ is a given constant.
\end{enumerate}
\end{definition}

\section{The main results}\label{main}

Our approach to solving the problem in Definition~\ref{Def1} will follow
the methodology for the analysis of stability and $L_2$-gain for networked
control systems proposed in~\cite{LFH-2012}. 

The proofs of the results are omitted due to space
limitation. The key technical tools used in those proofs are the Wirtinger's
inequality~\cite{IMA10} and the descriptor method~\cite{SCL01,SuplinAut07}.

As can be seen from (\ref{e.1}), if the
observer at node $i$ polls a channel at time $t_{k-p_i+1}$, the next time
the same channel will be polled at time $t_{k+1}$. The longest time between
polls of the same channel at node $i$ constitutes the maximum delay in
communication between node $i$ and its neighbours, which will be denoted
$\tau_i$:
\[
\tau_i=\max_{k} (t_{k+1}-t_{k-p_i+1}).
\]
The largest communication delay in the network is then
$\tau=\max_{i}\tau_i$. It is easy to see from these definitions that
$\tau=\max_{k} (t_{k+1}-t_{k-\bar p+1})$, where $\bar p=\max_i p_i$.

Consider the following Lyapunov-Krasovskii candidate for the system
(\ref{e.1}):
\begin{eqnarray}
  V_i(e_i)&=& e_i'Y_i^{-1}e_i+\int_{t-\tau_i}^t\!\! e^{-2\alpha_i(t-s)}
  e_i(s)'S_ie_i(s)ds
   \nonumber \\
   &+&\tau_i\int_{t-\tau_i}^t\!\! e^{-2\alpha_i(t-s)}\dot
   e_i(s)'(\tau_i+s-t)R_i\dot e_i(s)ds,\quad 
  \label{Vi}
\end{eqnarray}
where $Y_i=Y_i'>0$,
$R_i=R_i'\ge 0$, $S_i=S_i'\ge 0$ and
$\alpha_i\ge 0$, $i=1,\ldots,N$,
are matrices and constants to be
determined. $V_i(e_i)$ is a standard Lyapunov-Krasovskii functional used in the
literature on  exponential stability of systems with time-varying
delays; e.g., see \cite{LFH-2012}.

Given a matrix $W_i=W_i>0$, define 
\[
\mathcal{W}_i(u,z)={\pi^2\over 4}(u-z)'W_i(u-z).
\]

\begin{lemma}\label{dissip-ineq}
Suppose there exist gains $K_i$, $L_i$, matrices $W_i=W_i>0$, and constants
$\alpha_i>0$, $0<\pi_i<2\alpha_iq_i^{-1}$, $i=1,\ldots,N$, such that the
following vector dissipation inequality holds for all
$i=\ldots,N$: For $t\in[t_k,t_{k+1})$,
\begin{eqnarray}
  \lefteqn{\dot V_i(e_i)+2\alpha_iV_i(e_i)-\sum_{j\in \mathcal{V}_i} \pi_j
  V_j(e_j)} && \nonumber \\
&& +\left(\sum_{j:i\in \mathcal{V}_j}\tau_j^2\right)\dot e_i'W_i\dot e_i
 -\sum_{j\in \mathcal{V}_i}\mathcal{W}_j(e_j,e_j(t_{k-\nu_j^{k,i}+1}))
  \nonumber \\
&&  +\frac{1}{\gamma^2}(p_i+q_i)\|e_i\|^2 -\frac{2}{\gamma^2}e_i'\sum_{j\in
    \mathcal{V}_i} e_j - \|\xi_i\|^2\leq 0,
  \label{VLF}
\end{eqnarray}
where  and $\nu_j^{k,i}$ is the
index of $j$ in the ordered permutation set $\Pi^k(\mathcal{V}_i)$.
Then the system~(\ref{e.1}) satisfies conditions of Definition~\ref{Def1}.
\end{lemma}

In what follows we present a sufficient condition for the
dissipation inequality (\ref{VLF}) to hold.
We begin with a technical lemma which essentially restates the
corresponding lemma of~\cite{PKJ-2011} in the form convenient for the
subsequent use in the paper. Consider a vector
$\delta=[\delta_0',\ldots,\delta_{p_i}']'$, $\delta_\nu\in \mathbf{R}^n$.
Also, for given $n\times n$ matrices $R_i=R_i'\ge 0$ and $G_i$, define
\[
\Psi_i=\left[
  \begin{array}{cccc}
    R_i & \frac{1}{2}(G_i+G_i') & \ldots & \frac{1}{2}(G_i+G_i') \\
   \frac{1}{2}(G_i+G_i') & R_i  & \ldots & \frac{1}{2}(G_i+G_i') \\
    \vdots & \vdots & \ddots & \vdots \\
   \frac{1}{2}(G_i+G_i') & \frac{1}{2}(G_i+G_i')  & \ldots & R_i
 \end{array}\right].
\]

\begin{lemma}\label{Lem1}
Suppose the matrices $R_i=R_i'\ge 0$ and $G_i$ are such that
\begin{eqnarray}
\label{Park}
  \left[\begin{array}{cc}R_i & G_i\\ G_i' & R_i
    \end{array}\right]\ge 0.
\end{eqnarray}
Then
\begin{eqnarray*}
 \tau_i \left[
\frac{1}{t-t_k}\delta_0'R_i\delta_0
+ \sum_{\nu=1}^{p_i-1}\frac{1}{t_{k-\nu+1}-t_{k-\nu}}
\delta_\nu'R_i\delta_\nu\right.&&\\
\left.+\frac{1}{t_{k-p_i+1}-t+\tau_i }
\delta_{p_i}'R_i\delta_{p_i}\right]
&\ge& \delta' \Psi_i \delta.
\end{eqnarray*}
\end{lemma}

Next, we introduce a number of matrices. First, we introduce
\[
\bar\Psi_i=e^{-2\alpha_i\tau_i } T_i'\Psi_i T_i.
\]
it can be further partitioned  in accordance with the partition of
$\bar{\mathbf{e}}_i$: 
\[
\bar\Psi_i=\left[\begin{array}{ccc}
\bar\Psi_{i,11} & \bar\Psi_{i,12} & \bar\Psi_{i,13} \\
\bar\Psi_{i,12}' & \bar\Psi_{i,22} & \bar\Psi_{i,23} \\
\bar\Psi_{i,13}' & \bar\Psi_{i,23}' & \bar\Psi_{i,33} \\
\end{array}\right],
\]
Then we introduce the correspondingly partitioned matrix
\begin{eqnarray}
&&
\tilde\Psi_i=\left[\begin{array}{ccc}
\tilde\Psi_{i,11} & \tilde\Psi_{i,12} & \tilde\Psi_{i,13} \\
\tilde\Psi_{i,12}' & \tilde\Psi_{i,22} & \tilde\Psi_{i,23} \\
\tilde\Psi_{i,13}' & \tilde\Psi_{i,23}' & \tilde\Psi_{i,33} \\
\end{array}\right], \label{tildePsi}
\end{eqnarray}
where we let
\begin{eqnarray*}
\tilde\Psi_{i,11} &=& \bar\Psi_{i,11}-2\alpha_iY_i^{-1}-S_i, \\
\tilde\Psi_{i,33} &=& \bar\Psi_{i,33}+e^{-2\alpha_i\tau_i}S_i,\\
\tilde\Psi_{i,\mu\nu}&=&\bar\Psi_{i,\mu\nu} \quad \mbox{for all other
  elements of $\tilde\Psi_i$}.
\end{eqnarray*}
Also, the following matrices will be used in the sequel:
\begin{eqnarray*}
&&\bar \Phi_{i,11}= \left[\begin{array}{ccc}
\pi_{j_1}Y_{j_1}^{-1}+{\pi^2\over 4}W_{j_1} & \ldots & 0 \\
\vdots & \ddots & \vdots \\
0 & \ldots & \pi_{j_{p_i}}Y_{j_{p_i}}^{-1}+{\pi^2\over
4}W_{j_{p_i}}
\end{array} \right], \\
&&\bar \Phi_{i,22}=
\left[\begin{array}{cccc}
{\pi^2\over 4}W_{j_1}  & 0 & \ldots & 0 \\
0 & {\pi^2\over 4}W_{j_2}  & \ldots & 0 \\
\vdots & \vdots & \ddots & \vdots \\
0 & 0 & \ldots & {\pi^2\over 4}W_{j_{p_i}}
\end{array} \right], \\
&&\bar \Phi_{i,12}=\bar \Phi_{i,21}=-\bar \Phi_{i,22}.
\end{eqnarray*}

Finally, to formulate our first result concerned with the
analysis of consensus performance of the observer network~(\ref{UP7.C.d}),
we introduce the matrix
\begin{eqnarray}
\Xi_i=
\left[
\begin{array}{ccccccc}
\Xi_{aa} & \Xi_{ab} & \Xi_{ac} & 0 & 0 & \Xi_{af} & \Xi_{ag}\\
\star & \Xi_{bb} & \Xi_{bc} & -\tilde\Psi_{i,13} & \Xi_{be} & \Xi_{bf} & \Xi_{bg}\\
\star & \star & -\tilde\Psi_{i,22} & -\tilde\Psi_{i,23} & 0 & \Xi_{cf} & 0 \\
\star & \star & \star & -\tilde\Psi_{i,33} & 0 & 0 & 0 \\
\star & \star & \star & \star & -\bar\Phi_{i,11} &  -\bar\Phi_{i,12} & 0 \\
\star & \star & \star & \star & \star & \Xi_{ff} & \Xi_{fg}\\
\star & \star & \star & \star & \star & \star & -I
\end{array}\right], \nonumber \\
\label{Xi_i}
\end{eqnarray}
where we have used the following notation
\begin{eqnarray*}
\Xi_{aa}&=&\tau_i^2R_i+\left(\sum_{j:~i\in\mathcal{V}_j}\tau_j^2\right) W_i
-Z_i-Z_i', \\
\Xi_{ab}&=&Y_i^{-1}-X_i+Z_i'(A-L_iC_i), \\
\Xi_{ac}&=&-Z_i'(\mathbf{1}_{p_i}'\otimes K_iH_i),\\
\Xi_{af}&=&\mathbf{1}_{p_i}'\otimes (-Q_i+Z_i'K_iH_i),\\
\Xi_{ag}&=&Z_i'(B-L_iD_i),\\
\Xi_{bb}&=&
\frac{(p_i+q_i)}{\gamma^2}I-\tilde \Psi_{i,11} \\
&&+X_i'(A-L_iC_i)+(A-L_iC_i)'X_i,
\\
\Xi_{bc}&=&-\tilde\Psi_{i,12}-(\mathbf{1}_{p_i}'\otimes X_i'K_iH_i), \\
\Xi_{be}&=&-\frac{1}{\gamma^2}(\mathbf{1}_{p_i}'\otimes I), \\
\Xi_{bf}&=&\mathbf{1}_{p_i}'\otimes(X_i'K_iH_i+(A-L_iC_i)'Q_i),\\
\Xi_{bg}&=&X_i'(B-L_iD_i),\\
\Xi_{cf}&=&-\mathbf{1}_{p_i}\mathbf{1}_{p_i}'\otimes (H_i'K_i'Q_i),\\
\Xi_{ff}&=&\mathbf{1}_{p_i}\mathbf{1}_{p_i}'\otimes (Q_i'K_iH_i+H_i'K_i'Q_i)
-\bar \Phi_{i,22}, \\
\Xi_{fg}&=&\mathbf{1}_{p_i}\otimes Q_i'(B-L_iD_i).
\end{eqnarray*}
Here $X_i$, $Z_i$ and $Q_i$ are arbitrary $n\times n$ matrices. These
matrices are introduced in order to apply the descriptor
method~\cite{SCL01} to derive the following theorem.

\begin{theorem}\label{analysis}
Suppose there exist matrices $Y_i=Y_i'>0$, $X_i$, $Z_i$, $Q_i$,
$W_i=W_i'\ge 0$, $S_i=S_i'\ge 0$, $R_i=R_i'\ge 0$, $G_i$, constants
$\alpha_i>0$, $0\le \pi_i<2\alpha_i q_i^{-1}$, and gain matrices $K_i,L_i$,
$i=1,\ldots, N$, which satisfy  the LMI (\ref{Park}) and
\begin{eqnarray}
\label{LMI.analysis}
\Xi_i<0.
\end{eqnarray}
Then the corresponding observer network (\ref{UP7.C.d}) solves the problem posed in
Definition~\ref{Def1}. The
matrix $P$ in condition (\ref{objective.i.1}) 
corresponding to this solution is
$P=
\frac{1}{N}\sum_{i=1}^N(Y_i^{-1}
+S_i\frac{1-e^{-2\alpha_i\tau_i}}{2\alpha_i})$.    
\end{theorem}

Theorem~\ref{analysis} serves as the basis for derivation of 
the main result of this paper, given below in Theorem~\ref{T1}, which is 
a sufficient
condition for synthesis of distributed observer networks of the form
(\ref{UP7.C.d}). Consider the following matrix
\begin{eqnarray}
\bar\Xi_i=
\left[
\begin{array}{ccccccc}
\bar\Xi_{aa} & \bar\Xi_{ab} & \bar\Xi_{ac} & 0 & 0 & \bar\Xi_{af} & \bar\Xi_{ag}\\
\star & \bar\Xi_{bb} & \bar\Xi_{bc} & -\tilde\Psi_{i,13} & \bar\Xi_{be} & \bar\Xi_{bf} & \bar\Xi_{bg}\\
\star & \star & -\tilde\Psi_{i,22} & -\tilde\Psi_{i,23} & 0 & \bar\Xi_{cf} & 0 \\
\star & \star & \star & -\tilde\Psi_{i,33} & 0 & 0 & 0 \\
\star & \star & \star & \star & -\bar\Phi_{i,11} &  -\bar\Phi_{i,12} & 0 \\
\star & \star & \star & \star & \star & \bar\Xi_{ff} & \bar\Xi_{fg}\\
\star & \star & \star & \star & \star & \star & -I
\end{array}\right], \nonumber \\
\label{barXi_i}
\end{eqnarray}
where 
\begin{eqnarray*}
\bar\Xi_{aa}&=&\tau_i^2R_i+\left(\sum_{j:~i\in\mathcal{V}_j}\tau_j^2\right) W_i
-\epsilon_i X_i-\epsilon_i X_i', \\
\bar\Xi_{ab}&=&Y_i^{-1}-X_i+\epsilon_i (X_i'A-U_iC_i), \\
\bar\Xi_{ac}&=&-\epsilon_i (\mathbf{1}_{p_i}'\otimes F_iH_i),\\
\bar\Xi_{af}&=&\mathbf{1}_{p_i}'\otimes (-\bar\epsilon_i X_i+\epsilon_i
F_iH_i),\\
\bar\Xi_{ag}&=&\epsilon_i (X_i'B-U_iD_i),\\
\bar\Xi_{bb}&=&
\frac{(p_i+q_i)}{\gamma^2}I-\tilde \Psi_{i,11} \\
&&+X_i'A-U_iC_i+A'X_i-C_i'U_i',
\\
\bar\Xi_{bc}&=&-\tilde\Psi_{i,12}-\mathbf{1}_{p_i}'\otimes (F_iH_i), \\
\bar\Xi_{be}&=&-\frac{1}{\gamma^2}(\mathbf{1}_{p_i}'\otimes I), \\
\bar\Xi_{bf}&=&\mathbf{1}_{p_i}'
      \otimes(F_iH_i+\bar\epsilon_i AX_i-\bar\epsilon_i C'_iU_i'),\\
\bar\Xi_{bg}&=&X_i'B-U_iD_i),\\
\bar\Xi_{cf}&=&-\mathbf{1}_{p_i}\mathbf{1}_{p_i}'\otimes (\bar\epsilon_i H_i'F_i'),\\
\bar\Xi_{ff}&=&\bar\epsilon_i\mathbf{1}_{p_i}\mathbf{1}_{p_i}'
           \otimes (F_iH_i+H_i'F_i')
-\bar \Phi_{i,22}, \\
\bar\Xi_{fg}&=&\bar\epsilon_i\mathbf{1}_{p_i}\otimes (X_i'B-U_iD_i).
\end{eqnarray*}

\begin{theorem}\label{T1}
Suppose there exists matrices $Y_i=Y_i'>0$, $X_i$, $\det X_i\neq 0$,
$F_i$, $U_i$, $S_i=S_i'\ge 0$, $R_i=R_i'\ge 0$, $W_i=W_i'\ge 0$, $G_i$, and
constants $\alpha_i>0$, $0\ge \pi_i<2\alpha_iq_i^{-1}$, $\epsilon_i>0$,
$\bar\epsilon_i>0$, $i=1,\ldots, N$,
which satisfy the LMI (\ref{Park}) and
\begin{eqnarray}
\label{LMI.synthesis}
\bar \Xi_i<0.
\end{eqnarray}
Then the network of observers (\ref{UP7.C.d}) with
\begin{eqnarray}
K_i=(X_i')^{-1})F_i, \quad L_i=(X_i')^{-1}U_i,
\label{gains}
\end{eqnarray}
solves the distributed estimation problem posed in Definition~\ref{Def1}.
The matrix $P$ in condition (\ref{objective.i.1}) corresponding to this
solution is 
$P=
\frac{1}{N}\sum_{i=1}^N(Y_i^{-1}+S_i\frac{1-e^{-2\alpha_i\tau_i}}{2\alpha_i})$.
\end{theorem}

\section{Conclusions}\label{Concl}

The paper has presented a sufficient LMI condition for the design of a
Round-Robin 
interconnection protocol for networks of distributed observers. We have shown
that the proposed protocol allows one to use sampled-data communications
between the observers in the network, and does not require a combinatorial 
gain scheduling. As a result, the node observers are shown to be capable of
achieving the $H_\infty$ consensus objective introduced
in~\cite{U6,LaU1,U7}. 

\newcommand{\noopsort}[1]{} \newcommand{\printfirst}[2]{#1}
  \newcommand{\singleletter}[1]{#1} \newcommand{\switchargs}[2]{#2#1}

\end{document}